\documentclass{article}

\usepackage{amssymb,amsmath,epsfig,amscd,eucal,psfrag,amsthm}

\newtheorem{theorem}{Theorem}[section]
\newtheorem{lemma}[theorem]{Lemma}

\newtheorem{corollary}[theorem]{Corollary}
\newtheorem{remark}[theorem]{Remark}

\newtheorem{definition}[theorem]{Definition}

{\bf}{\rm}
\newtheorem{proposition}[theorem]{Proposition}

\newenvironment{pf}{\par\medskip\noindent{\em Proof. }}{\hfill $\square$\par\medskip}

\newcommand{\G}{\mathcal{G}}

\input xy
\xyoption{all}

\newtheorem{notation}[theorem]{Notation}

\def\PSL{\mbox{\rm{PSL}}}
\def\P{\mbox{\rm{P}}}

\def\SL{\mbox{\rm{SL}}}

\def\int{\mbox{\rm{int}}}

\def\PSL{\mbox{\rm{PSL}}}
\def\SL{\mbox{\rm{SL}}}

\def\r{{\mathbb R}}
\def\z{{\mathbb Z}}
\def\q{{\mathbb Q}}
\def\Fp{\mathbb{F}_p}
\def\c{{\mathbb C}}

\def\h3{{\mathbb H}^3}
\def\pp{{\mathfrak{p}}}
\def\qq{{\mathfrak{q}}}

\def\ok{{\cal{O}_{\rm{k}}}}
\def\okp{{\cal{O}_{{\rm{k}}_{\pp}}}}
\def\oLp{{\cal{O}_{{\rm{L}}_{\pp}}}}
\def\okq{{\cal{O}_{{\rm{k}}_{\qq}}}}

\title{Separability of double cosets and conjugacy classes in $3$-manifold groups}
\author{Emily Hamilton, Henry Wilton and Pavel Zalesskii}

\begin{document}

\maketitle

\begin{abstract}
Let $M = \h3 / \Gamma$ be a hyperbolic $3$-manifold of finite volume.  We show that if $H$ and $K$ are abelian subgroups of $\Gamma$ and $g \in \Gamma$, then the double coset $HgK$ is separable in $\Gamma$. As a consequence we prove that if M is a closed, orientable, Haken $3$-manifold and the fundamental group of every hyperbolic piece of the torus decomposition of $M$ is conjugacy separable then so is the fundamental group of $M$.  Invoking recent work of Agol and Wise, it follows that if $M$ is a compact, orientable 3-manifold then $\pi_1(M)$ is conjugacy separable.
\end{abstract}

\section{Introduction}
\label{intro}

The profinite topology on a group $\Gamma$ is the coarsest topology in which every homomorphism from
$\Gamma$ to a finite group is continuous.  When $\Gamma$ is the fundamental group of a manifold $M$,
the profinite topology on $\Gamma$ encodes the finite-sheeted covering spaces of $M$,
and as such is of great interest in the field of low-dimensional topology.

\begin{definition}
If a subset $X$ of $\Gamma$ is closed in the profinite topology then $X$ is called \emph{separable}.
Equivalently, for every $\gamma\in\Gamma - X$ there is a homomorphism $\phi$ from $\Gamma$
to a finite group such that $\phi(\gamma)\notin\phi(X)$.
\begin{enumerate}
\item  A group $\Gamma$ is \emph{residually finite} if the trivial subgroup is separable in $\Gamma$.
\item A group $\Gamma$ is \emph{conjugacy separable} if every conjugacy class in $\Gamma$ is separable.
\item A group $\Gamma$ is \emph{subgroup separable} or \emph{locally extended residually finite (LERF)}
if every finitely generated subgroup of $\Gamma$ is separable in $\Gamma$.
\item A group $\Gamma$ is \emph{double-coset separable} if for every pair $H, K$ of finitely generated subgroups
of $\Gamma$, and every $g \in \Gamma$ the double coset $HgK$ is separable.
\end{enumerate}
\end{definition}

Note that conjugacy separability and subgroup separability both imply residual finiteness,
while double-coset separability implies subgroup separability.

Hempel showed, using Thurston's Geometrization Theorem, that the fundamental groups of
Haken 3-manifolds are residually finite \cite{He}.  Scott proved that the fundamental groups
of surfaces and Seifert-fibered 3-manifolds are subgroup separable \cite{S},
and asked whether the same holds for all 3-manifold groups.  Burns, Karass and  Solitar \cite{burns_note_1987}
answered this question in the negative by giving an example of a graph manifold with non-subgroup-separable fundamental group,
but it follows from the main theorems of recent preprints of Agol \cite{agol_virtual_????} and Wise \cite{wise_preprint_????} that the fundamental groups of hyperbolic 3-manifolds are subgroup separable. In this paper,
we are concerned with conjugacy separability and double-coset separability in 3-manifold groups.

\subsection{Conjugacy separability}

Our first theorem extends the results of \cite{wilton_profinite_2010}, in which it was proved that the
fundamental groups of graph manifolds are conjugacy separable.

\begin{theorem}\label{t: Reducing cs to pieces of the torus decomposition}
Let $M$ be a closed, orientable, Haken 3-manifold  and let $N_1,\ldots,N_m$ be the pieces of the
torus decomposition of $M$.  If each $\pi_1(N_i)$ is conjugacy separable then $\pi_1(M)$ is conjugacy separable.
\end{theorem}

This is an important step in the proof that the fundamental group of any compact, orientable 3-manifold is conjugacy separable, as we now explain.  We follow the same broad strategy that Hempel used in his proof of residual finiteness,
although there are more difficult technical obstacles to overcome.

Let $M$ be a compact, orientable 3-manifold, possibly with boundary. We are interested in the question of  whether or not $\pi_1(M)$ is conjugacy separable.  If $M$ has boundary then, cutting along compressing discs and appealing to the fact that a free product of conjugacy separable groups is conjugacy separable \cite{stebe_residual_1970}, we may assume that $\partial M$ is incompressible.  Let $D$ be the double of $M$ along $\partial M$.  Then $\pi_1(M)$ is a retract of $\pi_1(D)$.  In particular, a pair of elements is conjugate in $\pi_1(M)$ if and only if it is conjugate in $\pi_1(D)$, and it follows that if $\pi_1(D)$ is conjugacy separable then so is $\pi_1(M)$. In this way, we can reduce to the case in which $M$ is orientable and closed.  Passing to the pieces of the Kneser--Milnor decomposition, and appealing again to the fact that a free product of conjugacy separable  groups is conjugacy separable, we can reduce further to the irreducible case.

The next step is to pass to the pieces of the torus decomposition of $M$, described by Jaco--Shalen \cite{jaco_seifert_1979} and Johannson \cite{johannson_homotopy_1979}.  This point is the heart of Hempel's argument.
He proves a gluing theorem that reduces the residual finiteness of $\pi_1(M)$ to the residual finiteness
of the fundamental groups of the pieces.   In the context of conjugacy separability,
Theorem \ref{t: Reducing cs to pieces of the torus decomposition} supplies the necessary gluing theorem.

This reduces the question of which 3-manifold groups are conjugacy separable to the geometric case.
If $M$ is a torus bundle over a circle then $\pi_1(M)$ is polycyclic, and so is conjugacy separable
by a theorem of Remeslennikov \cite{remeslennikov_conjugacy_1969}.  Martino proved that the fundamental
groups of Seifert-fibered $3$-manifolds are conjugacy separable \cite{martino_proof_2007}.

The hyperbolic case is much more difficult, but dramatic progress has been made recently.  The crucial concept is the notion of a \emph{special} group, introduced by Haglund and Wise \cite{haglund_special_2008}\footnote{By a \emph{special} group we shall mean a group that is the fundamental group of a compact A-special cube complex.  The reader is referred to \cite{haglund_special_2008} for definitions.}.   Minasyan proved that special groups are conjugacy separable \cite{minasyan_hereditary_2009}, while Chagas and the third author gave conditions under which conjugacy separability passes to finite extensions \cite{chagas_bianchi_2010}.   Taking these two results together, it follows that, if $N$ is a hyperbolic 3-manifold and $\pi_1(N)$ is \emph{virtually} special (that is, $\pi_1(N)$ has a special subgroup of finite index) then $\pi_1(N)$ is conjugacy separable.

Many hyperbolic 3-manifold and orbifold groups are known to be virtually special \cite{bergeron_virtual_2008,bergeron_hyperplane_????,chagas_bianchi_2010,chesebro_virtually_2009}, among them non-cocompact arithmetic and standard cocompact arithmetic lattices.  Wise has announced a proof that the fundamental group of any hyperbolic 3-manifold containing an embedded geometrically finite surface is virtually special \cite{wise_research_2009}; the heart of his proof is contained in \cite{wise_preprint_????}.  Very recently, Agol has given a proof that the fundamental group of any closed hyperbolic 3-manifold group is virtually special \cite{agol_virtual_????}.  His proof uses the work of \cite{wise_preprint_????}, as well as Kahn and Markovic's resolution of the Surface Subgroup Conjecture \cite{kahn_immersing_????} and an extension of the techniques of \cite{agol_residual_????}.

Combining the recent work of Wise and Agol with Theorem \ref{t: Reducing cs to pieces of the torus decomposition} we deduce that the fundamental group of every closed, orientable 3-manifold is conjugacy separable.  As described above, conjugacy separability for any compact, orientable 3-manifold, possibly with boundary, reduces to the closed case.  We therefore have a complete resolution of the question of conjugacy separability for the fundamental groups of compact, orientable 3-manifolds.

\begin{theorem}\label{thm: All 3-manifolds!}
If $M$ is any compact, orientable 3-manifold then $\pi_1(M)$ is conjugacy separable.
\end{theorem}

We emphasize again that Theorem \ref{thm: All 3-manifolds!} depends on the results of \cite{agol_virtual_????} and \cite{wise_preprint_????}.

\subsection{Double-coset separability}

In \cite{wilton_profinite_2010}, the second and third authors proved a combination theorem for conjugacy-separable groups,
and were able to check the hypotheses in the case when the vertex groups are the fundamental groups of Seifert-fibered 3-manifolds.
To prove Theorem \ref{t: Reducing cs to pieces of the torus decomposition}, we need to check the hypotheses of the combination theorem
for the fundamental groups of hyperbolic 3-manifolds of finite volume. In particular, we need to prove that double cosets of peripheral
subgroups of Kleinian groups of finite covolume are separable. In fact, we prove the following more general result.

\begin{theorem}[Theorem \ref{doublecoset}]\label{introdoublecoset}
Let $M = \h3 / \Gamma$ be a hyperbolic $3$-orbifold of finite volume. If $H$ and $K$ are abelian
subgroups of $\Gamma$ and $g \in \Gamma$, then the double coset $HgK = \{ hgk \ \mid h \in H, k \in K \}$ is separable in $\Gamma$.
\end{theorem}

In the closed case, this result can also be deduced from the subgroup separability of $\Gamma$ \cite{agol_virtual_????}.
Indeed, Minasyan proved that, if $G$ is a word-hyperbolic group (as $\Gamma$ is when $M$ is closed),
and every quasi-convex subgroup $H$ is separable in $G$, then for any pair of quasi-convex subgroups
$H$ and $K$ and any $g\in G$ the double coset $HgK$ is separable \cite{minasyan_separable_2006}.
(Abelian subgroups of a word-hyperbolic group are always quasi-convex.)
However, Minasyan's theorem has not been generalized to the case in which $M$ has boundary.

Haglund and Wise prove that double cosets of hyperplane subgroups are separable in virtually special groups---see \cite{haglund_special_2008},
to which the reader is also referred for the definition of a \emph{hyperplane subgroup}.
When $\Gamma$ is virtually special and $H$ is a hyperplane subgroup of $\Gamma$ then $\Gamma$ has a
finite-index subgroup that splits as an amalgamated product or HNN extension over $H$.
It follows that if $\Gamma$ is a Kleinian group of finite covolume then no hyperplane subgroup of $\Gamma$ is abelian.

\subsection{An outline}

This paper is structured as follows.  In Section \ref{algebra} we state the algebraic results
underlying the proof of Theorem \ref{introdoublecoset}.  In Section \ref{DoubleCoset} we prove
Theorem \ref{introdoublecoset} and discuss some implications.  In Section \ref{s: Conjugacy separability}
we recall the results of \cite{wilton_profinite_2010} and prove
Theorem \ref{t: Reducing cs to pieces of the torus decomposition}.

\subsection*{Acknowledgments}

The first author would like to thank Eric Brussel and Alan Reid
for useful conversations.  The second author would like to thank Stefan Friedl for some comments on an earlier draft; he is supported by an EPSRC Career Acceleration Fellowship. The third author is supported by CNPq.

\section{Algebraic preliminaries}
\label{algebra}

In this section we prove algebraic results that will be used
in the proof of Theorem~\ref{introdoublecoset}.  We assume standard terminology
of algebraic number theory.  For reference see \cite{J}.

\begin{notation}
By a {\it number field} we mean a finite
field extension of $\q$.  If $k$ is a number field,
let $\ok$ denote the ring of algebraic integers of $k$.
If $\pp$ is a non-zero prime ideal of $\ok$,
then we complete $k$ at $\pp$ to obtain
the local field $k_{\pp}$, with ring of
algebraic integers $\okp$.
The ring $\okp$
has a unique maximal ideal.
The quotient of $\okp$
by this maximal ideal is called the {\it residue class field} of
$\okp$.  The quotient map
is called the {\it residue class field map}
with respect to $\pp$.
\end{notation}

We begin by stating  
two theorems and two corollaries from \cite{H3}.

\begin{theorem}
\label{theoremA}
Let $k$ be a number
field and let $\delta$ be a non-zero element of $k$
that is not a root of unity.  Let $S$ be a finite
set of prime ideals of $\ok$.
Then there exists a positive integer $n$ with
the following property.  For each integer $m \geq n$,
there exists a non-zero
prime ideal $\pp$ of $\ok$, lying outside of $S$,
such that $\delta \in \okp$ and the multiplicative
order of the image of $\delta$ in the residue class
field of $\okp$ is equal to $m$.
\end{theorem}

For the proof of Theorem~\ref{theoremA} please see Theorem 2.3
of \cite{H3}.

\begin{corollary}
\label{orderm}
Let $R$ be a finitely generated ring
in a number field $k$, let $\delta$ be
a non-zero element of $R$
that is not a root of unity, and let
$x_1, x_2, \ \ldots \ , x_j$ be non-zero
elements of $R$.
Then there exists a positive integer $n$
with the following property.  For each integer
$m \geq n$, there exist a finite field $F$
and a ring homomorphism $\eta: R \rightarrow F$
such that the multiplicative order of $\eta(\delta)$ is
equal to $m$ and $\eta(x_i) \neq 0$, for each
$1 \leq i \leq j$.
\end{corollary}

The deduction of Corollary~\ref{orderm} from Theorem~\ref{theoremA}
is given in \cite{H3}. (See Corollary 2.5 of \cite{H3}.)  We include the proof
in this paper for the convenience of the reader.

\begin{pf}
Fix a finite generating set $G$ of $R$.
Let $S$ denote the finite set of
prime ideals of $\ok$
which divide an element of $\{ G,
x_1, x_2, \ ... \ , x_j \}$.
By Theorem~\ref{theoremA}, there exists a positive
integer $n$ with the following property.
For each integer $m \geq n$ there exists a non-zero
prime ideal $\pp$ of $\ok$, lying outside of $S$, such that
$\delta \in \okp$ and the multiplicative
order of the image of $\delta$ in the residue class field
of $\okp$ is equal to $m$.  Fix $m \geq n$ and let
$\pp \subset \ok$ be the corresponding prime ideal.
Let $F$ denote the residue class field of
$\okp$ and let $\eta: \okp \rightarrow F$ denote
the residue class field map with respect to $\pp$.
Since $\pp \notin S$, $R \in \okp$ and
$x_1, x_2, \ ... , x_j$ are units in $\okp$.
Therefore, the restriction of $\eta$ to $R$ satisfies
the conclusion of the corollary.
\end{pf}

\begin{theorem}
\label{multringsep}
Let $k$ be a number field.  Let $\lambda$ and
$\omega$ be non-zero elements of $k$ such that
$\lambda$ is not a multiplicative
power of $\omega$.
Let $P$ be a finite set of prime ideals of $\ok$.
Then there exist primes $\pp$ and $\qq$,
lying outside of $P$,
such that $\lambda, \omega \in \okp
\cap \okq$ and $(\eta_{\pp} \times
\eta_{\qq})(\lambda)$ is not a multiplicative power of
$(\eta_{\pp} \times \eta_{\qq})(\omega)$.
\end{theorem}

For the proof of Theorem~\ref{multringsep} please see Theorem 2.7
of \cite{H3}.

\begin{corollary}
\label{notpower}
Let $R$ be a finitely generated
ring in a number field $k$.
Let $\lambda$ and
$\omega$ be non-zero elements of $R$ such that
$\lambda$ is not a multiplicative
power of $\omega$.  Then there exist a finite
ring $S$ and a ring homomorphism
$\eta: R \rightarrow S$ such that
$\eta(\lambda)$ is not a multiplicative power
of $\eta(\omega)$.
\end{corollary}

The deduction of Corollary~\ref{notpower} from Theorem~\ref{multringsep} is similar to the deduction of Corollary~\ref{orderm} from Theorem~\ref{theoremA}. See Corollary 2.8 of \cite{H3} for details.

Theorem~\ref{multringsep}
can be interpreted as a multiplicative
subgroup separability result for finitely generated rings lying in number fields.
We conclude this section with an additive
subgroup separability result.

\begin{theorem}
\label{addringsep}
Let $R$ be a finitely generated ring in a number field $k$.
By fixing a $\q$ embedding of $k$ into $\c$, we may view $k \subset \c$.
Let $\beta$ be an element of $R$ and set
$A = \{ m + n\beta \mid m, n \in \z \}$ and $B = \{ m + n\beta \mid m, n \in \q \}$.
If $b \in R - B$, then there exist a finite ring $S$
and a ring homomorphism $\eta: R \rightarrow S$ such that
$\eta(b) \notin \eta(A)$.
\end{theorem}

\begin{pf}
We first consider the case where $b \notin \q(\beta)$.
Let $L$ denote the normal closure of $k$ over $\q$.
Since $b \notin \q(\beta)$, there exists an element
$\sigma \in \mathrm{Gal}(L/\q)$ such that
$\sigma(b) \neq b$ and $\sigma$ fixes $\q(\beta)$ pointwise.
By the Tchebotarev Density Theorem,
there are infinitely many primes $p$ of $\q$ with unramified
extension $\pp$ in $L$ such that $\sigma$
is the Frobenius automorphism for ${\pp} / p$.
Fix one such $\pp / p$ such that
$R \subset \oLp$, where
$\oLp$ denotes the ring of integers
in the $\pp$-adic field $L_{\pp}$.
Such a choice is possible since $R$ is finitely
generated, and each generator is an integer in
$L_{\pp}$ for all but finitely many $\pp$.
Let $F_{\pp}$ denote the residue class field of
$\oLp$ and let
$\Fp$ denote the finite field of $p$ elements.
Let $\eta$ be the composition of the inclusion map of
$R$ into $\oLp$ with the residue map:
$$ \eta: R \hookrightarrow \oLp \rightarrow F_{\pp}.$$
Since $\sigma$ is the Frobenius automorphism of $L/\q$
with respect to $\pp / p$, $\mathrm{Gal}(L_{\pp} / \q_{p})
= \langle \sigma^{\prime} \rangle$ where
$\sigma^{\prime} = \sigma$ on $L$.  Since $A \subset
\q(\beta)$ and $\sigma$ fixes $\q(\beta)$ pointwise,
$A \subset \q_{p}$.  Since $\sigma(b) \neq b$,
$b \notin \q_{p}$.
The Galois group of $F_{\pp} / \Fp$ is also induced by $\sigma$.
It follows that $\eta(A) \subset \Fp$, but $\eta(b) \notin \Fp$.
Therefore, $\eta$ and $S = \Fp$ satisfy the conclusion of the theorem.

Now consider the case where $b \in \q(\beta)$.
Let $f$ be the minimal monic polynomial of $\beta$ over $\q$
and let $n$ be the degree of $f$.  Our assumption that
$b \notin B$ implies that $n > 2$.  Moreover, we can
express $b = a_0 + a_1\beta + a_2 \beta^2 + \ \ldots \ + a_{n-1}\beta^{n-1}$,
where $a_i \in \q$ and at least one coefficient $a_{i_0} \in
\{a_2, a_3, \ \ldots \ , a_{n-1} \}$ is non-zero.
Since $R$ is a finitely generated ring consisting of
algebraic numbers, $R$ is integral over $\z[1/s]$, for all but
finitely many integers $s$.  Choose $s$ such that
the coefficients of $f$ are in $\z[1/s]$ and $a_i \in \z[1/s],
\ \forall i \in \{ 0, 1, 2, \ \ldots \ , n-1 \}$.
Let $\z[1/s][T]$ denote the polynomial ring of $\z[1/s]$.
Let $p$ be a prime that does not divide $s$ or the numerator
of $a_{i_0}$ and let $\Fp$ denote the finite field of $p$ elements.
Since $f$ is a monic polynomial, the map
$$\z[1/s][T] \rightarrow \z[1/s][\beta], \ {\mathrm{given \ by}}
\ T \rightarrow \beta,$$ is an epimorphism with kernel $(f)$.  Therefore, the map
$$\rho: \z[1/s][\beta] \rightarrow \z[1/s][T]/(f), \ {\mathrm{given \ by}}
\ \beta \rightarrow T + (f),$$ is an isomorphism.
The quotient map $$\z[1/s] \rightarrow \z[1/s]/(p) \cong \Fp$$
induces
$$\phi: \z[1/s][T]/(f) \rightarrow {\Fp}[T]/(\overline{f}),$$
where $\overline{f}$ is the image of $f$ in the polynomial ring ${\Fp}[T]$.
Write $\overline{f} = \overline{f}_1 \overline{f}_2 \ \ldots \ \overline{f}_m$
as a product of irreducible factors in $\Fp[T]$, and let
$F_i = {\Fp}[T]/(\overline{f}_i)$. Then the natural map
$$\psi: \Fp[T]/(\overline{f}) \rightarrow F_1 \times F_1 \times \ \ldots \ \times F_m$$
is an isomorphism. Let
$$\rho_i: \z[1/s][\beta] \rightarrow \z[1/s][T]/(f) \rightarrow
{\Fp}[T]/(\overline{f}) \rightarrow F_1 \times F_1 \times \ \ldots \
\times F_m \rightarrow F_i$$ be the composition of $(\psi \phi
\rho)$ with the projection of $F_1 \times F_1 \times \ \ldots \
\times F_m$ onto $F_i$. Since $F_i$ is a field, the kernel of
$\rho_i$ is a maximal ideal $\pp_i$ of $\z[1/s][\beta]$.  By
replacing $R$ with $R[1/s]$, if necessary, we may assume that $1/s
\in R$.  Thus, $R$ is an integral extension of $\z[1/s][\beta]$.
Therefore, there exists a maximal ideal $\qq_i$ of $R$ such that
$(\z[1/s][\beta]) \cap \qq_i = \pp_i$.  Let $\eta_i$ denote the
quotient map $R \rightarrow R/{\qq_i}$. Note that the restriction
of $\eta_i$ to $\z[1/s][\beta]$ is equal to $\rho_i$. Let $\pi_i:
R/{\qq_1} \times R/{\qq_2} \times \ \ldots \ \times R/{\qq_m}
\rightarrow R/{\qq_i}$ denote the projection onto the $i$-th
factor. By the universal property of direct products of rings,
there exists a unique ring homomorphism $\eta: R \rightarrow
R/{\qq_1} \times R/{\qq_2} \times \ \ldots \ \times R/{\qq_m}$ such
that $\pi_i \circ \eta = \eta_i$. We claim that $\eta(b) \notin
\eta(A)$. Assume, to the contrary, that $\eta(b) = \eta(c_0 +
c_1\beta$) for some $c_0 + c_1\beta \in A$.  Then
\[
b - (c_0 + c_1\beta) =  (a_0 - c_0) + (a_1 - c_1)\beta + a_2 \beta^2 + \ \ldots \ + a_{n-1}\beta^{n-1} \]
is in the kernel of $\rho_i$ for every $i
\in \{1, 2, \ \ldots \ , m \}$. Let
\[
h(T) = (a_0 - c_0) + (a_1 - c_1)T + a_2 T^2 + \ \ldots \ + a_{n-1} T^{n-1} \in \z[1/s][T]
\]
and let	 $\overline{h}$ be the image of $h$ in ${\Fp}[T]$. Then
$\overline{f}$ divides $\overline{h}$ in $\Fp[T]$. But this
contradicts the fact that the image of $a_{i_0}$ is non-zero in
$\Fp$ and so $1 < {\mathrm{deg}}(\overline{h}) < n$. Therefore,
$\eta$ and $S = R/{\qq_1} \times R/{\qq_2} \times \ \ldots \ \times
R/{\qq_m}$ satisfy the conclusion of the theorem.
\end{pf}

\section{Proof of double-coset separability}
\label{DoubleCoset}

In this section we prove that double cosets of abelian subgroups
of Kleinian groups of finite covolume are separable.  The proof will
use the following proposition from \cite{N}.

\begin{proposition}
\label{niblo}
Let $G_0, H$ and $K$ be finitely generated subgroups of a group $G$
and set $H_0 = H \cap G_0$ and $K_0 = K \cap G_0$.  If $[G: G_0]$
is finite and if $H_0K_0$ is separable in $G_0$, then $HK$ is separable in $G$.
\end{proposition}

For a proof of this result see Proposition 2.2 of \cite{N}.

\begin{theorem}
\label{doublecoset}
Let $M = \h3 / \Gamma$ be an orientable hyperbolic $3$-orbifold
of finite volume.  Let $H$ and $K$ be abelian subgroups of $\Gamma$,
and let $g \in \Gamma$.
Then the double coset $HgK = \{ hgk \ \mid h \in H, k \in K \}$
is separable in $\Gamma$.
\end{theorem}

\begin{pf}
As noted in \cite{N}, since the profinite topology on $\Gamma$
is equivariant under left and right multiplication, the double
coset $HgK$ is closed in $\Gamma$ if and only if $H^gK = g^{-1}HgK$
is closed in $\Gamma$.  Therefore, to prove the theorem, it suffices to
show that if $H$ and $K$ are abelian subgroups of $\Gamma$,
then the double coset $HK$ is separable in $\Gamma$.
Note that $g \in HK$ if and only if $g^{-1} \in KH$.  Therefore,
$HK$ is separable in $\Gamma$ if and only if $KH$ is separable
in $\Gamma$.

By Selberg's Lemma \cite{A}, $\Gamma$ has a subgroup of
finite index $\Gamma_0$ which is torsion free.
Let $H_0 = H \cap \Gamma_0$ and $K_0 = K \cap \Gamma_0$.
By Proposition~\ref{niblo}, if $H_0K_0$ is separable in
$\Gamma_0$, then $HK$ is separable in $\Gamma$.  Therefore,
we may assume that $\Gamma$ is torsion free, and thus the
fundamental group of a hyperbolic manifold.
Since abelian subgroups of finitely generated Kleinian groups
are separable \cite{AH}, we may assume that $H$ and $K$
do not commute.  In particular, both
$H$ and $K$ are non-trivial. Since $\Gamma$ is torsion free,
the non-trivial abelian subgroups of $\Gamma$ are free abelian
of rank $1$ or $2$.  The free abelian subgroups of rank $1$
can be generated by loxodromic or parabolic isometries.  The
free abelian subgroups of rank $2$ are generated
by parabolic isometries.  If $H$ has rank $1$, let $H^{\prime}$
be the maximal cyclic subgroup of $\Gamma$ containing $H$.  If
$H$ has rank $2$, let $H^{\prime}$ be the maximal abelian subgroup
of $\Gamma$ containing $H$.  Then $H$ has finite index in $H^{\prime}$.
Let $\{ a_1, a_2, \ \ldots \ , a_{n-1} \}$ be a set of non-trivial coset
representatives of $H^{\prime} / H$.
By \cite{AH}, $H$ is a separable subgroup of $\Gamma$.
Therefore, there exists a subgroup $\Gamma_H$ of finite index in $\Gamma$
such that $H \subset \Gamma_H$ but $\Gamma_H \cap
\{ a_1, a_2, \ \ldots \ , a_{n-1} \} = \emptyset$.
Then $\Gamma_H \cap H^{\prime} = H$.  In a similar way,
define $K^{\prime}$ and choose $\Gamma_K$ of finite index in $\Gamma$
such that $\Gamma_K \cap K^{\prime} = K$.  Let $\Gamma_0 = \Gamma_H \cap
\Gamma_K$ and set $H_0 = H \cap \Gamma_0$ and $K_0 = K \cap \Gamma_0$.
Then $\Gamma_0$ is a subgroup of finite index in $\Gamma$,
$H_0 = \Gamma_0 \cap H^{\prime}$ and $K_0 = \Gamma_0 \cap K^{\prime}$.
By Proposition~\ref{niblo}, by replacing $\Gamma$
with $\Gamma_0$, if necessary, we may assume that rank-one elements
of $\{ H, K \}$ are maximal cyclic subgroups of $\Gamma$ and rank-two
elements of $\{ H, K \}$ are maximal abelian subgroups of $\Gamma$.
This assumption will be used if $H$ and/or $K$ is parabolic.  However, the
proof does not require loxodromic subgroups to be maximal.

Given the assumptions and reductions above, we need to consider the following cases.

\bigskip

\noindent
{\it Case 1.}  $H$ loxodromic, $K$ loxodromic

\medskip

The group of orientation preserving isometries
of $\h3$ may be identified with $\PSL(2,\c)$.
Thus there exists a discrete, faithful representation
$\rho: \pi_1(M) \rightarrow \PSL(2,\c)$
which is well-defined up to conjugation in $\PSL(2,\c)$.
Let $\q (\mathrm{tr}\Gamma)$ denote the field
obtained by adjoining the traces of the elements of
$\rho(\Gamma)$ to $\q$. Since $M$ has finite volume, it follows
from Mostow Rigidity that $\q (\mathrm{tr}\Gamma)$ is a
number field.  By Proposition 2.2(e) of \cite{B} we may
conjugate $\rho(\Gamma)$ in $\PSL(2,\c)$ to lie in a finite field
extension of $\q (\mathrm{tr}\Gamma)$.  Therefore we
may view $\Gamma \subset \PSL(2,k)$, where $k$ is a finite
extension of $\q$.
Write $H = \langle h \rangle$ and $K = \langle k \rangle$.
In this case, $h$ and $k$ are conjugate
in $\PSL(2,\c)$ to matrices of the form
$$\pm \begin{pmatrix}
\lambda & 0 \\
0 & \lambda^{-1} \\
\end{pmatrix},
\ \lambda \in \c, \ \vert \lambda \vert \neq 1, \
\mathrm{and} \
\pm \begin{pmatrix}
\omega & 0 \\
0 & \omega^{-1} \\
\end{pmatrix},
\ \omega \in \c, \ \vert \omega \vert \neq 1, $$
respectively.  Since $\q (\mathrm{tr}\Gamma)$ is a number field,
the eigenvalues, $\lambda$ and $\omega$, are algebraic numbers.
By adjoining $\lambda$ and $\omega$ to $k$, if necessary,
we may assume that $h$ and $k$ are diagonalizable over $k$.
Therefore, after conjugating in $\PSL(2,k)$,
$$h = \pm \begin{pmatrix}
\lambda & 0 \\
0 & \lambda^{-1} \\ \end{pmatrix} \
\mathrm{and} \
k = \pm \begin{pmatrix}
a & b \\ c & d \\ \end{pmatrix}
\begin{pmatrix}
\omega & 0 \\ 0 & \omega^{-1} \\ \end{pmatrix}
{\begin{pmatrix}
a & b \\ c & d \\ \end{pmatrix}}^{-1},$$
for some $$g = \pm \begin{pmatrix}
a & b \\ c & d \\ \end{pmatrix} \in \PSL(2,k).$$
By assumption, $h$ and $k$
do not commute. Since $\Gamma$ is discrete,
it follows that the fixed points of $h$
on the sphere at infinity are disjoint from the fixed points of $k$.
Therefore, all of the elements in $\{ a, b, c, d \}$ are non-zero.

Let $G$ be the subgroup of $\PSL(2,k)$
generated by $\Gamma$ and $g$.  We will show that
$HKg$ is separable in $G$.
Since the profinite topology on $G$
is equivariant under left and right multiplication,
it will follow that $HK$ is separable in $G$, implying
that $HK$ is separable in $\Gamma$, as required.
To prove that $HKg$ is separable in $G$,
note that
$$HKg = \Big\lbrace \pm \begin{pmatrix}
\lambda^m & 0 \\
0 & \lambda^{-m} \\ \end{pmatrix}
\begin{pmatrix}
a & b \\ c & d \\ \end{pmatrix}
\begin{pmatrix}
\omega^n & 0 \\ 0 & \omega^{-n} \\ \end{pmatrix}$$
$$= \pm \begin{pmatrix}
a\lambda^m\omega^n & b \lambda^m\omega^{-n} \\
c\lambda^{-m}\omega^n & d\lambda^{-m}\omega^{-n} \\
\end{pmatrix}  \ \Big\vert \ m, n \in \z \Big\rbrace,$$
and let
$$\gamma = \pm \begin{pmatrix}
r & s \\ t & u \\ \end{pmatrix} \in G
- HKg$$ be given.
Since $M$ has finite volume, $\Gamma$ is finitely generated.
Therefore, $G$ is finitely generated.
Let $R$ be the ring generated by the coefficients of
$G$ over $\z$.  Then $G \subset
\PSL(2,R) \subset \PSL(2,k)$.
Suppose that rs/ab is not a multiplicative power
of $\lambda^2$.  Then by Corollary~\ref{notpower},
there exist a finite ring $S$ and a ring homomorphism
$\eta: R \rightarrow S$ such that $\eta(ab) \neq 0$ and $\eta$(rs/ab)
is not a multiplicative power of $\eta(\lambda^2)$.
The map $\eta$ induces a group homomorphism
$$\overline{\eta}:G
\hookrightarrow \PSL(2,R) \rightarrow \PSL(2,S).$$
Suppose that $\overline{\eta}(\gamma) \in
\overline{\eta}(HKg)$.  Then there exist elements $m, n \in \z$ such that
$\overline{\eta}(\gamma) = \overline{\eta}(h^mk^ng).$
Equating coefficients, we have:
$$\Big\lbrace \begin{array}{clcr}
\eta(r) = \eta(a\lambda^m\omega^n),
& \eta(s) = \eta(b\lambda^m\omega^{-n}), \cr
\eta(t) = \eta(c\lambda^{-m}\omega^n),
& \eta(u) = \eta(d\lambda^{-m}\omega^{-n})
\cr \end{array} \Big\rbrace $$
\centerline{or}
$$\Big\lbrace \begin{array}{clcr}
\eta(-r) = \eta(a\lambda^m\omega^n),
& \eta(-s) = \eta(b\lambda^m\omega^{-n}), \cr
\eta(-t) = \eta(c\lambda^{-m}\omega^n),
& \eta(-u) = \eta(d\lambda^{-m}\omega^{-n})
\cr \end{array} \Big\rbrace.$$
In either case, $\eta($rs) $= \eta$(ab$\lambda^{2m})$, a contradiction.
Therefore, we may
assume that rs/ab $= \lambda^{2m_0}$, for some $m_0 \in \z$.
By a similar argument, we may assume that
rt/ac $= \omega^{2n_0}$, for some $n_0 \in \z$.

By assumption, $\gamma \notin HKg$.  In particular,
$$\begin{pmatrix}
r & s \\ t & u \\ \end{pmatrix} \notin
\Big\lbrace \pm \begin{pmatrix}
a\lambda^{m_0}\omega^{n_0} & b \lambda^{m_0}\omega^{-n_0} \\
c\lambda^{-m_0}\omega^{n_0} & d\lambda^{-m_0}\omega^{-n_0} \\
\end{pmatrix} \Big\rbrace.$$   Since $R$ is finitely generated,
$R \subset \okp$ for all but finitely many prime ideals
$\pp$ of $\ok$.  For each of these primes $\pp$ the residue
class field map $\eta_{\pp}: \okp \rightarrow F_{\pp}$ induces
a group homomorphism $$\overline{\eta}_{\pp}:G
\hookrightarrow \PSL(2,\okp) \rightarrow \PSL(2,F_{\pp}),$$
where $F_{\pp}$ is the residue class field of $\okp$.
Choose $\pp$ such that
$$ (\ast) \ \ \begin{pmatrix}
\eta_{\pp}(r) & \eta_{\pp}(s) \\ \eta_{\pp}(t) & \eta_{\pp}(u) \\ \end{pmatrix}
\notin \Big\lbrace \pm \begin{pmatrix}
\eta_{\pp}(a\lambda^{m_0}\omega^{n_0}) & \eta_{\pp}(b\lambda^{m_0}\omega^{-n_0}) \\
\eta_{\pp}(c\lambda^{-m_0}\omega^{n_0}) & \eta_{\pp}(d\lambda^{-m_0}\omega^{-n_0}) \\
\end{pmatrix} \Big\rbrace.$$  Suppose $\overline{\eta}_{\pp}(\gamma) \in
\overline{\eta}_{\pp}(HKg)$.  Then there exist elements $m, n \in \z$
such that $\overline{\eta}_{\pp}(\gamma) = \overline{\eta}_{\pp}(h^mk^ng).$
Equating coefficients as above,
$\eta_{\pp}(\lambda^{2m}) = \eta_{\pp}$(rs/ab)
$= \eta_{\pp}(\lambda^{2m_0})$
and $\eta_{\pp}(\omega^{2n}) = \eta_{\pp}$(rt/ac)
$= \eta_{\pp}(\omega^{2n_0})$.
Therefore, $\eta_{\pp}(\lambda^m) = \pm \eta_{\pp}(\lambda^{m_0})$
and $\eta_{\pp}(\omega^n) = \pm \eta_{\pp}(\omega^{n_0})$,
contradicting $(\ast)$.  This completes the proof in this case.

\bigskip

\noindent
{\it Case 2.}  $H$ loxodromic, $K$ maximal parabolic

\medskip

The proof of this case follows from the proof
of Lemma 3.2 in \cite{H2}.

\bigskip

\noindent
{\it Case 3.}  $H$ maximal parabolic, $K$ maximal parabolic

\medskip

Since $M$ is a hyperbolic $3$-manifold of finite volume,
we may view $M$ as the interior of a compact manifold $M^{\prime}$  with
a finite number of tori boundary components.
Such a manifold can be obtained from $M$ by truncating the cusp tori.
Furthermore, $\pi_1(M^{\prime}) \cong \pi_1(M) = \Gamma$.
There is a one-to-one correspondence between the boundary tori of $M^{\prime}$
and the conjugacy classes of maximal parabolic subgroups of $\Gamma$.
Suppose that $H$ and $K$ correspond to the same
boundary component of $M^{\prime}$.  Then there exists an  element $\zeta \in \Gamma$
such that $K = \zeta H \zeta^{-1}$.  Since $H$ and $K$ do
not commute, $\zeta \notin H$.
Since $H$ is a separable subgroup
of $\Gamma$ \cite{LN}, there exist a finite group $G$ and a group homomorphism
$f: \Gamma \rightarrow G$ such that $f(\zeta) \notin f(H)$.  Let $\Gamma_0$
be kernel of $f$, and set $H_0 = H \cap \Gamma_0$ and $K_0 = K \cap \Gamma_0$.
Since $\Gamma_0$ is normal in $\Gamma$, $K_0 = \zeta H_0 \zeta^{-1}$.
If there exists an element $\upsilon \in \Gamma_0$ such that
$K_0 = \upsilon H_0 \upsilon^{-1}$, then $(\upsilon^{-1}\zeta) H_0
(\upsilon^{-1}\zeta)^{-1} = H_0$. Therefore, $\upsilon^{-1}\zeta$
fixes the parabolic fixed point of $H$.  Since $H$ is a maximal
parabolic subgroup of $\Gamma$, $\upsilon^{-1}\zeta \in H$.  But then $f(\zeta)
\in f(H)$, a contradiction.  We conclude that $H_0$ and $K_0$ are not
conjugate in $\Gamma_0$.  By Proposition~\ref{niblo},
if $H_0K_0$ is separable in
$\Gamma_0$, then $HK$ is separable in $\Gamma$.  Therefore, by replacing
$\Gamma$ with $\Gamma_0$, if necessary, we may assume that $H$
and $K$ correspond to different boundary components of $M^{\prime}$.
Let $T_1$ and $T_2$ denote the boundary components
of $M^{\prime}$ corresponding to $H$ and $K$, respectively.
By Thurston's Hyperbolic Dehn Surgery
Theorem \cite{BP}, we may choose generators
$h_1, h_2$ for $H$  such that
$M_{h_1}$ and $M_{h_2}$ are complete
hyperbolic $3$-manifolds of finite volume,
where $M_{h_1}$ and $M_{h_2}$ are the manifolds obtained by Dehn
surgery on $M^{\prime}$ along $T_1$ sending $h_1$ and $h_2$, respectively,
to a meridian of the attached solid torus.  Similarly, we may choose
generators $k_1, k_2$ for $K$ such that $M_{k_1}$ and $M_{k_2}$
are complete hyperbolic $3$-manifolds of finite volume, where
$M_{k_1}$ and $M_{k_2}$ are obtained by Dehn surgery on $M^{\prime}$ along $T_2$.

For simplicity we use the fact that
the representation $\rho: \Gamma \rightarrow \PSL(2,\c)$
may be lifted to a representation
$$\Gamma \rightarrow \SL(2,\c).$$
(See Proposition 3.1.1 of \cite{CS}.)  Therefore, we view
$\Gamma \subset \SL(2,\c)$.
By an argument similar to that in Case 1,
we may assume that $\Gamma \subset
\SL(2,k)$, where $k$ is a finite field extension of $\q$,
$$ h_1 = \begin{pmatrix} 1 & 1 \\ 0 & 1 \\ \end{pmatrix} \ \mathrm{and}
\ h_2 = \begin{pmatrix} 1 & \beta_1 \\ 0 & 1 \\ \end{pmatrix},$$
for a fixed element $\beta_1 \in \c - \r$.  If $A_1 = \{ m + n\beta_1
\ \vert \ m, n \in \z \}$, then
$$H = \Big\{ \begin{pmatrix} 1 & x \\ 0 & 1 \\ \end{pmatrix}
\ \Big\vert \ \ x \in A_1  \Big\}.$$
Moreover, $$k_1 = \begin{pmatrix} a_1 & b_1 \\ c_1 & d_1 \\ \end{pmatrix}
\begin{pmatrix} 1 & 1 \\ 0 & 1 \\ \end{pmatrix}
\begin{pmatrix} a_1 & b_1 \\ c_1 & d_1 \\ \end{pmatrix}^{-1}
\ \mathrm{and} \
k_2 = \begin{pmatrix} a_1 & b_1 \\ c_1 & d_1 \\ \end{pmatrix}
\begin{pmatrix} 1 & \beta_2 \\ 0 & 1 \\ \end{pmatrix}
\begin{pmatrix} a_1 & b_1 \\ c_1 & d_1 \\ \end{pmatrix}^{-1},$$
for fixed elements $$\begin{pmatrix} a_1 & b_1 \\ c_1 & d_1 \\ \end{pmatrix}
\in \SL(2, \c) \ \mathrm{and} \ \beta_2 \in \c - \r.$$
If $A_2 = \{ m + n \beta_2 \ \vert \ m, n \in \z \}$, then
$$ K = \Big\{ \begin{pmatrix} 1 - ya_1c_1 & y{a_1}^2 \\ -y{c_1}^2 &
1 + ya_1c_1 \\ \end{pmatrix}
\ \Big\vert \ \ y \in A_2 \Big\}.$$
Since $H$ and $K$ do not commute, $c_1 \neq 0$.
Note that ${a_1}^2$ and ${c_1}^2$ are elements
of the coefficient field $k$.  Since $k$ is a number field, $a_1$
and $c_1$ are algebraic numbers.  Therefore, after adjoining $a_1$
and $c_1$ to $k$, and conjugating $\Gamma$ in $\SL(2,k)$
by $$\begin{pmatrix} 1 & a_1/c_1 \\ 0 & 1 \\ \end{pmatrix},$$
we may assume that
$$ H =
\Big\{ \begin{pmatrix} 1 & -a_1/c_1 \\ 0 & 1 \\ \end{pmatrix}
\begin{pmatrix} 1 & x \\ 0 & 1 \\ \end{pmatrix}
\begin{pmatrix} 1 & a_1/c_1 \\ 0 & 1 \\ \end{pmatrix} =
\begin{pmatrix} 1 & x \\ 0 & 1 \\ \end{pmatrix} \
\Big\vert \ \ x \in A_1 \Big\}, \ \mathrm{and}$$
$$ K = \Big\{ \begin{pmatrix} 1 & -a_1/c_1 \\ 0 & 1 \\ \end{pmatrix}
\begin{pmatrix} 1 - ya_1c_1 & y{a_1}^2 \\ -y{c_1}^2 & 1 + ya_1c_1 \\ \end{pmatrix}
\begin{pmatrix} 1 & a_1/c_1 \\ 0 & 1 \\ \end{pmatrix} =
\begin{pmatrix} 1 & 0 \\ -y{c_1}^2 & 1 \\ \end{pmatrix} \
\Big\vert \ \ y \in A_2 \Big\}.$$
Then $$ HK = \Big\{
\begin{pmatrix} 1 & x \\ 0 & 1 \\ \end{pmatrix}
\begin{pmatrix} 1 & 0 \\ -y{c_1}^2 & 1 \\ \end{pmatrix} =
\begin{pmatrix} 1 - xy{c_1}^2 & x \\ -y{c_1}^2 & 1 \\ \end{pmatrix}
\Big\vert \ \ x \in A_1, \ y \in A_2 \Big\}.$$

To show that $HK$ is separable in $\Gamma$,
let $$\gamma = \begin{pmatrix}
r & s \\ t & u \\ \end{pmatrix} \in \Gamma
- HK$$ be given.  Suppose that $u \neq 1$.
As in Case 1, there exists a finitely generated ring $R_1$
such that $\Gamma \subset \SL(2,R_1) \subset \SL(2,k)$.  Since
$R_1$ is finitely generated, $R_1 \subset \okp$ for all but finitely
many primes $\pp$ of $\okp$.
For each of these primes $\pp$ the residue
class field map $\eta_{\pp}: \okp \rightarrow F_{\pp}$ induces
a group homomorphism $$\overline{\eta}_{\pp}:\Gamma
\hookrightarrow \SL(2,\okp) \rightarrow \SL(2,F_{\pp}),$$
where $F_{\pp}$ is the residue class field of $\okp$.
Choose $\pp$ such that $\eta_{\pp}(u) \neq \eta_{\pp}(1)$.  Then
$\overline{\eta}_{\pp}(\gamma) \notin
\overline{\eta}_{\pp}(HK)$, as required. Therefore, we
may assume that $u = 1$.  Since the determinant of $\gamma$
is equal to $1$,
$$ \gamma = \begin{pmatrix} 1 + st & s \\ t & 1 \\ \end{pmatrix}.$$
Let $B_1 = \{ m+n\beta_1 \ \vert \ m, n \in \q \}$
and $B_2 = \{ m+n\beta_2 \ \vert \ m, n \in \q \}$.
Suppose that either $s \notin B_1$ or $-t/{c_1}^2 \notin B_2$.
By Theorem~\ref{addringsep}, there exist
a finite ring $S$ and a ring homomorphism $\eta: R_1 \rightarrow S$
such that $\eta(s) \notin \eta(B_1)$ or $\eta(-t/{c_1}^2) \notin \eta(B_2)$,
respectively.  This ring homomorphism induces a group homomorphism
$$\overline{\eta}:\Gamma
\hookrightarrow \SL(2,R_1) \rightarrow \SL(2,S),$$
such that $\overline{\eta}(\gamma) \notin
\overline{\eta}(HK)$.  Therefore, we may assume that
$s \in B_1$ and $-t/{c_1}^2 \in B_2$.  Moreover, since
$\gamma \notin HK$, either $s \notin A_1$ or $-t/{c_1}^2 \notin A_2$.
We will assume that $s \notin A_1$.  (The argument if $-t/{c_1}^2 \notin A_2$
is similar, with the roles of $H$ and $K$ interchanged).
Since $s \in B_1$, there exists a non-zero integer $v_0$ such
that $v_0s \in A_1$.  Write $v_0s = m_0 + n_0 \beta_1$, where $m_0, n_0 \in \z$.
Since $s \notin A_1$, either $v_0$ does not divide $m_0$ or $v_0$
does not divide $n_0$.
We will assume that $v_0$ does not divide $m_0$.
(The argument if $v_0$ does not divide $n_0$ is similar, with the roles
of $h_1$ and $h_2$ interchanged).
Let $M_{h_2}$ be the hyperbolic $3$-manifold obtained from $M^{\prime}$
by Dehn surgery along $T_1$, as defined above.
Let $$\phi: \Gamma \cong \pi_1(M^{\prime}) \rightarrow \pi_1(M_{h_2})$$
be the homomorphism induced by inclusion.  Then
$\phi(H)$ is an infinite cyclic loxodromic subgroup of $\pi_1(M_{h_2})$
generated by $\phi(h_1)$, and $\phi(h_2)$ is trivial.
By assumption, $H$ and $K$ correspond to different boundary components
of $M^{\prime}$.  Therefore,
$\phi(K)$ is a maximal parabolic subgroup of $\pi_1(M_{h_2})$.
If $\phi(\gamma) \notin \phi(HK)$, then we are done by Case 2.
Therefore, we may assume that
$\phi(\gamma) = \phi(h_1^{m_1})\phi(k_{\gamma})$, for some
$m_1 \in \z$ and $k_{\gamma} \in K$.
Since $v_0$ does not divide $m_0, \ v_0m_1-m_0 \neq 0$.

Recall that $\Gamma \subset \SL(2,R_1) \subset \SL(2, k)$.
Let $L$ denote the normal closure of $k$ over $\q$ and let $\tau
\in \mathrm{Gal}(L/\q)$ represent complex conjugation.
Since $\beta_1 \in \c - \r$, $\tau(\beta_1) \neq \beta_1$.
By the Tchebotarev Density Theorem,
there exist infinitely many primes $p$ of $\q$ with unramified
extension $\pp$ in $L$ such that $\tau$
is the Frobenius automorphism for ${\pp} / p$.
Fix one such $\pp / p$ such that $p$ is an odd prime,
$p$ does not divide $ v_0m_1 - m_0$ and
$R_1 \subset \oLp$, where
$\oLp$ denotes the ring of integers
in the $\pp$-adic field $L_{\pp}$.
Let $F_{\pp}$ denote the residue class field of
$\oLp$ and let
$\Fp$ denote the finite field of $p$ elements.
Let $\eta_1$ be the composition of the inclusion map of
$R_1$ into $\oLp$ with the residue map:
$$ \eta_1: R_1 \hookrightarrow \oLp \rightarrow F_{\pp}.$$
Since $\tau$ is the Frobenius automorphism of $L/\q$
with respect to $\pp / p$, $\mathrm{Gal}(L_{\pp} / \q_{p})
= \langle \tau^{\prime} \rangle$ where
$\tau^{\prime} = \tau$ on $L$.
Since $\tau(\beta_1) \neq \beta_1$,
$\beta_1 \notin \q_{p}$.
The Galois group of $F_{\pp} / \Fp$ is also induced by $\tau$.
It follows that $\eta_1(\beta_1) \notin \Fp$.
The map $\eta_1: R_1 \rightarrow F_{\pp}$ induces
a group homomorphism $$\psi_1: \Gamma
\hookrightarrow \SL(2,R_1) \rightarrow \SL(2,F_{\pp}).$$
If $\psi_1(\gamma) \notin \psi_1(HK)$, then we are done.
Suppose that $\psi_1(\gamma) = \psi_1(hk)$ for some
$h \in H, \ k \in K$.  If $h = h_1^m h_2^n$
and $x = m + n\beta_1$,
then $$\begin{pmatrix}
\eta_1(1 + st) & \eta_1(s) \\
\eta_1(t) & 1 \\ \end{pmatrix} =
\psi_1(\gamma) = \psi_1(hk) =
\begin{pmatrix}
\eta_1(1 - xy{c_1}^2) & \eta_1(x) \\
\eta_1(-y{c_1}^2) & 1 \\ \end{pmatrix},$$
for some $y \in A_2$ corresponding to $k$.
Therefore, $\eta_1(v_0m + v_0n\beta_1) = \eta_1(v_0x) = \eta_1(v_0s) = \eta_1(m_0 + n_0 \beta_1)$.
Since $\eta_1(\beta_1) \notin \Fp$, the set $\{ 1, \eta_1(\beta_1) \}$
is linearly independent over $\Fp$.
It follows that $$ (\ast) \ \  v_0m \equiv m_0 \ (\mathrm{mod} \ p).$$

Now consider the hyperbolic $3$-manifold $M_{h_2}$.
Recall that $\phi(H)$ is a loxodromic subgroup of $\pi_1(M_{h_2})$
generated by $\phi(h_1)$,  $\phi(K)$ is a maximal parabolic
subgroup of $\pi_1(M_{h_2})$,
and $\phi(\gamma) = \phi(h_1^{m_1} k_{\gamma})$, for some $m_1 \in \z$
and $k_{\gamma} \in K$.  As before, there exists a finitely generated ring $R_2$ in
a number field $F$, such that $\pi_1(M_{h_2}) \subset \SL(2,R_2)
\subset \SL(2,F)$.  Moreover, after conjugating $\pi_1(M_{h_2})$
in $\SL(2,F)$, if necessary, we may assume that
$$\phi(h_1) =  \begin{pmatrix}
\lambda & 0 \\
0 & \lambda^{-1} \\ \end{pmatrix} \
\mathrm{and}$$
$$\phi(K) =  \Big\{ \begin{pmatrix}
a_2 & b_2 \\ c_2 & d_2 \\ \end{pmatrix}
\begin{pmatrix}
1 & z \\ 0 & 1 \\ \end{pmatrix}
{\begin{pmatrix}
a_2 & b_2 \\ c_2 & d_2 \\ \end{pmatrix}}^{-1} \ \Big\vert \
z \in  A_3 = \{ m + n \beta_3 \ \vert \ m, n \in \z \} \Big\},$$
for some $$  \begin{pmatrix}
a_2 & b_2 \\ c_2 & d_2 \\ \end{pmatrix} \in \SL(2,F) \ \mathrm{and} \
\beta_3 \in \c - \r.$$
Then $$\phi(HK) = \Big\{ \begin{pmatrix} \lambda^m (1 - a_2c_2z) &
\lambda^m {a_2}^2z \\
-\lambda^{-m}{c_2}^2z &
\lambda^{-m}(1 + a_2{c_2}z) \\ \end{pmatrix} \
\Big\vert \ m \in \z, \ z \in A_3 \Big\}.$$
Write $$\phi(k_{\gamma}) =
\begin{pmatrix} 1 & z_{\gamma} \\ 0 & 1 \\ \end{pmatrix}, \ z_{\gamma}
\in A_3.$$

Since $\pi_1(M_{h_2})$ is discrete, $\phi(H)$ and $\phi(K)$
do not share a fixed point.  Therefore, $a_2$ and $c_2$ are non-zero.
By Corollary~\ref{orderm}, there exist a finite field $S$
and a ring homomorphism $\eta_2: R_2 \rightarrow S$
such that the multiplicative order of $\eta_2(\lambda)$
is divisible by $p$, and $\eta_2(a_2)$ and $\eta_2(c_2)$
are non-zero in $S$.
Since $\lambda$ is a unit in $R_2$, $\eta_2(\lambda)
\neq 0$.  Let $o$ denote the multiplicative
order of $\eta_2(\lambda)$.  Consider the composition:
$$\psi_2:  \Gamma \rightarrow \pi_1(M_{h_2}) \hookrightarrow
\SL(2,R_2) \rightarrow \SL(2,S),$$
where the first map  $\phi: \Gamma \rightarrow \pi_1(M_{h_2})$
is induced by inclusion and the last map $\SL(2,R_2) \rightarrow \SL(2,S)$ is
induced by $\eta_2$.
If $\psi_2(\gamma) \notin \psi_2(HK)$, then we are done.
Suppose that $\psi_2(\gamma) = \psi_2(hk)$, for some
$h \in H, \ k \in K$.  If $h = h_1^m h_2^n$, then
$$\psi_2(\gamma) =
 \begin{pmatrix} \eta_2(\lambda^{m_1} (1 - a_2c_2z_{\gamma})) &
\eta_2(\lambda^{m_1} {a_2}^2z_{\gamma}) \\
\eta_2(-\lambda^{-{m_1}}{c_2}^2z_{\gamma}) &
\eta_2(\lambda^{-{m_1}}(1 + a_2{c_2}z_{\gamma})) \\ \end{pmatrix} =$$
$$ \psi_2(hk) = \begin{pmatrix} \eta_2(\lambda^m (1 - a_2c_2z)) &
\eta_2(\lambda^m {a_2}^2z) \\
\eta_2(-\lambda^{-m}{c_2}^2z) &
\eta_2(\lambda^{-m}(1 + a_2{c_2}z)) \\ \end{pmatrix},$$
for some $a \in A_3$ corresponding to $k$.
This gives the equations:
$$ \begin{array}{clcr}
\eta_2(\lambda^{m_1} (1 - a_2c_2z_{\gamma})) =
\eta_2(\lambda^m (1 - a_2c_2z)) \cr
\eta_2(\lambda^{m_1} {a_2}^2z_{\gamma}) = \eta_2(\lambda^m {a_2}^2z) \cr
\eta_2(-\lambda^{-{m_1}}{c_2}^2z_{\gamma}) =
\eta_2(-\lambda^{-m}{c_2}^2z)  \cr
\eta_2(\lambda^{-{m_1}}(1 + a_2{c_2}z_{\gamma})) =
\eta_2(\lambda^{-m}(1 + a_2{c_2}z)).
\cr \end{array}$$  If $\eta_2(z_{\gamma}) = 0$, then $\eta_2(z) = 0$,
and so $\eta_2(\lambda^m) = \eta_2(\lambda^{m_1})$.  If
$\eta_2(z_{\gamma}) \neq 0$, then by solving for $\eta_2(z/z_{\gamma})$
in the second and third equations, we have that
$\eta_2(\lambda^{2m}) = \eta_2(\lambda^{2m_1})$.  In either case,
$2m \equiv 2m_1 \ (\mathrm{mod} \ o)$.
Since $p$ divides $o$ and $p$ is an odd prime, it follows that
$$(\ast \ast) \ \ m \equiv m_1 \ (\mathrm{mod} \ p).$$

Finally, consider the product
$$(\psi_1 \times \psi_2): \Gamma \rightarrow \SL(2, F_{\pp}) \times
\SL(2, S).$$  If $(\psi_1 \times \psi_2)(\gamma) \notin
(\psi_1 \times \psi_2)(HK)$, then we are done.  Suppose that
$(\psi_1 \times \psi_2)(\gamma) = (\psi_1 \times \psi_2)(hk)$
for some $h \in H, k \in K$.  Then $\psi_1(\gamma) =
\psi_1(hk)$ and $\psi_2(\gamma) = \psi_2(hk)$.  If $h = h_1^m h_2^n$,
then by $(\ast)$ and $(\ast \ast)$,
$v_0m \equiv m_0 \ (\mathrm{mod} \ p)$ and
$m \equiv m_1 \ (\mathrm{mod} \ p).$
Therefore, $v_0m_1 \equiv m_0 \ (\mathrm{mod} \ p)$.
This contradicts the fact that we chose $p$ not to
divide $v_0m_1 - v_0$.  This completes the proof of Case 3.

\bigskip

\noindent
{\it Case 4.}  $H$ loxodromic, $K$ parabolic

\medskip

Let $K^{\prime}$ be the maximal parabolic subgroup of $\Gamma$
containing $K$.  By Case 2 and the assumptions at the beginning of the proof,
we may assume that $K$ is a maximal cyclic subgroup of $\Gamma$.  Therefore,
there exists a basis $ \{ k_1, k_2 \}$ of $K^{\prime}$ such that
$K = \langle k_1 \rangle$.  Let $H = \langle h \rangle $.
To prove that $HK$ is separable in $\Gamma$,
let $\gamma \in \Gamma - HK$ be given.  By Case 2, we may assume
that $\gamma \in HK^{\prime} - HK$.  Therefore, $\gamma = h^{a_0} k_1^{m_0}
k_2^{n_0}$, for some $a_0, m_0, n_0 \in \z$.  Since
$\gamma \notin HK$, $n_0 \neq 0$.
As above, we may assume that
$\Gamma \subset \SL(2,R) \subset \SL(2, k)$, where $R$
is a finitely generated ring contained in a number field $k$.
Moreover, $$ k_1 = \begin{pmatrix} 1 & 1 \\ 0 & 1 \\ \end{pmatrix} \
\mathrm{and} \ k_2 = \begin{pmatrix} 1 & \beta \\ 0 & 1 \\
\end{pmatrix},$$ for some $\beta \in \c - \r$.
Since $H$ is loxodromic, $h$ is conjugate in $\SL(2, \c)$
to an element
$$ f = \begin{pmatrix} \lambda & 0 \\ 0 & {\lambda}^{-1} \\ \end{pmatrix},
\ \lambda \in \c, \ \vert \lambda \vert \neq 1.$$
The eigenvalue $\lambda$ is an algebraic number.
By adjoining $\lambda$ to $k$, we may assume that $h$ is diagonalizable
over $\SL(2,k)$.  Therefore, $$ h =
\begin{pmatrix} a & b \\ c & d \\ \end{pmatrix}
\begin{pmatrix} \lambda & 0 \\ 0 & \lambda^{-1} \\ \end{pmatrix}
{\begin{pmatrix} a & b \\ c & d \\ \end{pmatrix}}^{-1}, \
\mathrm{for \ some} \ \begin{pmatrix} a & b \\ c & d \\ \end{pmatrix}
\in \SL(2,k) .$$ By expanding $R$, if necessary, we may assume
that $\{ \lambda, {\lambda}^{-1}, a, b, c, d \} \subset R$.
Let $L$ denote the normal closure of $k$ over $\q$ and let $\tau
\in \mathrm{Gal}(L/\q)$ represent complex conjugation.
As in Case 3, by the Tchebotarev Density Theorem,
there exist infinitely many primes $p$ of $\q$ with unramified
extension $\pp$ in $L$ such that $\tau$
is the Frobenius automorphism for ${\pp} / p$.
Fix one such $\pp / p$ such that
$p$ does not divide $n_0$
and $R \subset \oLp$, where
$\oLp$ denotes the ring of integers
in the $\pp$-adic field $L_{\pp}$.
Let $F_{\pp}$ denote the residue class field of
$\oLp$ and let
$\Fp$ denote the finite field of $p$ elements.
Let $\eta$ be the composition of the inclusion map of
$R$ into $\oLp$ with the residue map:
$$ \eta: R \hookrightarrow \oLp \rightarrow F_{\pp}.$$
Since $\tau$ is the Frobenius automorphism of $L/\q$
with respect to $\pp / p$, $\eta(\beta) \notin \Fp$.
The map $\eta: R \rightarrow F_{\pp}$ induces
a group homomorphism $$\psi: \Gamma
\hookrightarrow \SL(2,R) \rightarrow \SL(2,F_{\pp}).$$
Suppose that $\psi(\gamma) \in \psi(HK)$.  Then
$\psi(h^{a_0} k_1^{m_0} k_2^{n_0}) = \psi(h^a k_1^m)$,
for some $a, m \in \z$. Therefore, $\psi(k_1^{m_0 - m}k_2^{n_0})
= \psi(h^{a - a_0})$. The trace of each element in $K^{\prime}$
is equal to $2$.
It follows that $\eta(\lambda^{a - a_0}) + 1/\eta(\lambda^{a - a_0}) = 2$,
and so $\eta(\lambda^{a - a_0}) = 1$.
This means that $\psi(f^{a - a_0})$ is trivial.  Since $h$ is conjugate to
$f$ in $\SL(2, R)$, $\psi(h^{a - a_0}) =
\psi(k_1^{m_0 - m} k_2^{n_0})$ is trivial.
Therefore, $\eta(m_0 - m + n_0 \beta) = 0$.  Since $\eta(\beta) \notin \Fp$,
the set $\{ 1, \eta(\beta) \}$ is linearly independent over $\Fp$.
It follows that $\eta(n_0) = 0$.  But this contradicts the fact
that $p$ does not divide $n_0$.  We conclude that
$\psi(\gamma) \notin \psi(HK)$, as required

\bigskip

\noindent
{\it Case 5.}  $H$ parabolic, $K$ parabolic

\medskip

Let $H^{\prime}$ and $K^{\prime}$ be the maximal parabolic subgroups
of $\Gamma$, containing $H$ and $K$, respectively.
By the assumptions at the beginning of the proof,
either $H = H^{\prime}$, or $H$ is a maximal
cyclic subgroup of $H^{\prime}$.  A similar statement is true for $K$.
If $H = H^{\prime}$ and $K = K^{\prime}$ then we are done by Case 3.
Therefore, we assume that either $H$ or $K$ is infinite
cyclic. Choose bases $\{ h_1, h_2 \}$ for $H^{\prime}$ and $\{ k_1, k_2 \}$
for $K^{\prime}$ such that $H = \langle h_1 \rangle$ if $H$
is infinite cyclic, and $K = \langle k_1 \rangle$ if $K$ is infinite cyclic.
To show that $HK$ is separable in $\Gamma$, let
$\gamma \in \Gamma - HK$ be given.   By Case 3, we may assume that
$\gamma \in H^{\prime} K^{\prime} - HK$.  Therefore, $\gamma =
h_1^{m_1} h_2^{n_1} k_1^{m_2} k_2^{n_2}$, for some $m_1, m_2, n_1, n_2
\in \z$.  Since $\gamma \notin HK$, either ({\it i}) $H$ is infinite cyclic
and $n_1 \neq 0$ or ({\it ii}) $K$ is infinite cyclic and $n_2 \neq 0$.
Without loss of generality, assume that ({\it i}) holds.

The argument is then very similar to the argument in Case 4.  We
may assume that $\Gamma \subset \SL(2, R) \subset \SL(2,k)$, where
$R$ is a finitely generated ring in a number field $k$.
Moreover, $$ h_1 = \begin{pmatrix} 1 & 1 \\ 0 & 1 \\ \end{pmatrix},
\ h_2 = \begin{pmatrix} 1 & \beta_1 \\ 0 & 1 \\
\end{pmatrix}, \ \mathrm{and}$$
$$K = \Big\{ \begin{pmatrix} a & b \\ c & d \\ \end{pmatrix}
\begin{pmatrix} 1 & y \\ 0 & 1 \\ \end{pmatrix}
{\begin{pmatrix} a & b \\ c & d \\ \end{pmatrix}}^{-1}
= \begin{pmatrix} 1 - acy & a^2 y \\
-c^2 y & 1 + acy \\ \end{pmatrix} \ \Big\vert \
y \in \{ m + n \beta_2 \} \Big\},$$
$$\mathrm{for \ some} \ \begin{pmatrix} a & b \\ c & d \\ \end{pmatrix}
\in \SL(2,R) \ \mathrm{and} \ \beta_1, \beta_2 \in \c - \r.$$
By assumption, $H$ and $K$ do not commute.  Therefore $c \neq 0$.
As in Case 4, by the Tchebotarev Density Theorem, there exist
a prime $p$ that does not divide $n_1$,
a finite field $F_{\pp}$ of characteristic $p$, and a ring homomorphism
$$\eta: R \rightarrow F_{\pp}$$ such that the set $\{ 1, \eta(\beta_1) \}$
is linearly independent over $\Fp$,
and $\eta(c) \neq 0$.  The map $\eta$ induces a group
homomorphism $$\psi: \Gamma \hookrightarrow \SL(2,R) \rightarrow \SL(2, F_{\pp}).$$
Suppose that $\psi(\gamma) \in \psi(HK)$.  Then
$\psi(h_1^{m_1} h_2^{n_1} k_1^{m_2} k_2^{n_2}) =
\psi(h_1^{u_1} k_1^{u_2} k_2^{v_2}),$
for some $u_1, u_2, v_2 \in \z$. Therefore,
$\psi(h_1^{m_1 - u_1} h_2^{n_1}) =
\psi(k_1^{u_2 - m_2} k_2^{v_2 - n_2})$.
Since $H$ is upper triangular and $\eta(c) \neq 0$,
it follows that $\psi(k_1^{u_2 - m_2} k_2^{v_2 - n_2})
= \psi(h_1^{m_1 - u_1} h_2^{n_1})$ is trivial.
Therefore, $\eta(m_1 - u_1 + n_1\beta_1) = 0$.
Since $\{ 1, \eta(\beta_1) \}$ is linearly independent
over $\Fp$, $\eta(n_1) = 0$.  But this contradicts the fact
that $p$ does not divide $n_1$.
\end{pf}

We conclude this section with some corollaries to
the proof of Theorem~\ref{doublecoset}.
Since only Case 3 of Theorem~\ref{doublecoset}
uses the full strength of the finite covolume assumption,
it is natural to consider finitely generated Kleinian groups
which are not necessarily of finite covolume.

\begin{corollary}
Let $\Gamma$ be a finitely generated, torsion-free Kleinian group.
Given an abelian subgroup $G$ of $\Gamma$, let $A(G)$ denote the
maximal abelian subgroup of $\Gamma$ containing $G$.
Suppose that $H$ and $K$ are abelian subgroups of $\Gamma$
such that $A(H)$ and $A(K)$ are not both free abelian of rank $2$.
Then the double coset $HK$ is separable in $\Gamma$.
\end{corollary}

\begin{pf}
If $\Gamma$ is elementary, then $\Gamma$ is virtually
abelian and hence the result follows from Proposition~\ref{niblo}
and \cite{AH}.  If $\Gamma$ is non-elementary, then
$\Gamma$ is isomorphic to a geometrically finite
Kleinian group $\Gamma_1$ such that
({\it i}) every maximal parabolic
subgroup of $\Gamma_1$ has rank $2$, and ({\it ii})
the traces of the elements of $\Gamma_1$ are algebraic
numbers. (See Theorem 1 of \cite{AH} or Theorem 4.2 of \cite{Sh}.)
By replacing $\Gamma$ with $\Gamma_1$, we assume that conditions
({\it i}) and ({\it ii}) hold for $\Gamma$.
Our assumptions then imply that at least one
of $H$ or $K$ must be loxodromic.  The proof then follows
from Cases 1, 2 and 4 of Theorem~\ref{doublecoset}.
\end{pf}

In Theorem~\ref{doublecoset}, we prove that certain double cosets
of Kleinian groups of finite covolume are closed in the profinite
topology on $\Gamma$.  We now consider
the congruence topology on $\Gamma$.

\begin{definition}
Let $k$ be a number field and let
$\Gamma$ be a finitely generated subgroup of $\PSL(2,k)$.
Since $\Gamma$ is finitely generated, $\Gamma \subset \PSL(2,R) \subset \PSL(2,k)$,
where $R$ is a ring obtained from  $\ok$ by inverting a finite
number of elements.  This ring $R$ is Dedekind and, therefore, for
any non-zero ideal $I$, the quotient $R/I$ is finite.
The quotient map $R \rightarrow R/I$ induces a \emph{congruence homomorphism}
$$\eta: \Gamma \hookrightarrow \PSL(2,R) \rightarrow \PSL(2,R/I).$$
We say that a subset $X$ of $\Gamma$ is \emph {closed in the congruence topology}
on $\Gamma$ if for every element $\gamma \in \Gamma - X$,
there exists a congruence homomorphism $\eta$ such that
$\eta(\gamma) \notin \eta(X)$.
\end{definition}

Recall that a set $X \subset \Gamma$ is closed
in the profinite topology on $\Gamma$ if for every
element $\gamma \in \Gamma - X$, there exist a finite
group $G$ and a group homomorphism $\phi: \Gamma \rightarrow G$
such that $\phi(\gamma) \notin \phi(X)$. For the
profinite topology we consider all group homomorphisms from
$\Gamma$ into finite groups.  For the congruence
topology we consider only congruence homomorphisms.  Therefore, the
congruence topology is weaker than the profinite topology.

\begin{corollary}
\label{congruencetop}
Let $M = \h3 / \Gamma$ be an orientable hyperbolic $3$-orbifold
of finite volume.  Let $\rho: \Gamma \rightarrow  \PSL(2,\c)$
be a discrete, faithful representation such that the coefficient
field of $\rho(\Gamma)$ is a finite field extension of $\q$.
Then we may view $\rho(\Gamma) \subset \PSL(2, R)$,
where $R$ is a finitely generated ring in a number field $k$.
Let $H$ and $K$ be abelian subgroups of $\Gamma$,
and let $g \in \Gamma$. The double coset $HgK$ is closed in the congruence
topology on $\Gamma$, with respect to $\rho$ and $R$,
if one of the following conditions is satisfied.
\begin{itemize}
\item The groups $H$ and $K$ are both loxodromic subgroups of $\Gamma$.
\item Exactly one of $\{ H, K \}$ is a loxodromic subgroup of
$\Gamma$ and the other is a maximal parabolic subgroup of $\Gamma$.
\item Exactly one of $\{ H, K \}$ is a loxodromic subgroup of
$\Gamma$ and the other is a maximal cyclic parabolic subgroup of
$\Gamma$.
\end{itemize}
\end{corollary}

\begin{remark}
In Corollary~\ref{congruencetop}, we do not require any of the loxodromic
subgroups to be maximal.
\end{remark}

\begin{pf}
Since $M = \h3 / \Gamma$ is a Kleinian group of finite covolume,
there exists a discrete, faithful representation from $\Gamma$ into $\PSL(2, \c)$
such that the coefficient field of the image of $\Gamma$ is a finite field extension of $\q$.
By fixing one such representation $\rho$ we may view $\Gamma \subset \PSL(2, R)
\subset \PSL(2,k)$, where $R$ is a finitely generated ring in a number field $k$.
In the proof of Theorem~\ref{doublecoset}, given an abelian subgroup
$H$ of $\Gamma$, we adjoin finitely many algebraic numbers to $R$, if
necessary, and then conjugate $\Gamma$ in $\PSL(2,R)$ such that $H$
has a nice form.   We need to justify that we can replace
our original representation $\rho$ with the new representation
given by conjugation. (It is fine to replace our ring $R$
with a larger ring since any ring homomorphism from
the larger ring restricts to a ring homomorphism from $R$.)
To see this let $\alpha \in \PSL(2, R)$ and consider $\Gamma^{\prime}
= \alpha \Gamma \alpha^{-1} \subset \PSL(2,R)$.  Given a subset $X \subset \Gamma$
and an element $\gamma \in \Gamma - X$, let 
$X^{\prime} = \alpha X \alpha^{-1} \subset \Gamma^{\prime}$
and $\gamma^{\prime} = \alpha \gamma \alpha^{-1} \in
\Gamma^{\prime} - X^{\prime}$.  Suppose there exist
a finite ring $S$ and a ring homomorphism $R \rightarrow S$,
such that, under the induced group homomorphism $$\eta: \PSL(2, R)
\rightarrow \PSL(2, S),$$ $\eta(\gamma^{\prime}) \notin \eta(X^{\prime})$.
Then restricting $\eta$ to $\Gamma$ gives a congruence
homomorphism
$$ \Gamma \hookrightarrow \PSL(2, R) \rightarrow \PSL(2, S)$$
such that $\eta(\gamma) \notin \eta(X)$.
We conclude that, for the purposes of our proof, it is legitimate
to adjoin finitely many algebraic numbers to $R$
and then replace $\Gamma$ with a conjugate of $\Gamma$ in $\PSL(2,R)$.

We first show that loxodromic subgroups, maximal parabolic
subgroups, and maximal cyclic parabolic subgroups
are closed in the congruence topology on $\Gamma$,
with respect to $\rho$ and $R$.
This is well-known and the proof of some of the cases is in \cite{AH}.
We include the proof here for the convenience of the reader.

Let $A$ be a maximal abelian subgroup of $\Gamma$.
As above, we view $\Gamma \subset \PSL(2,R) \subset \PSL(2,k)$,
where $R$ is a finitely generated ring in a number field $k$.
After adjoining elements to $R$ and $k$, if necessary, we
may conjugate $\Gamma$ in $\PSL(2,R)$ such that $A$
is upper triangular.  Let
$$ \gamma =  \pm \begin{pmatrix}
a & b \\ c & d \\ \end{pmatrix}
\ \in \Gamma - A$$ be given.  Since $\Gamma$ is discrete and
$A$ is a maximal abelian subgroup of $\Gamma$, $\gamma$ and $A$ do not
share a fixed point on the sphere at infinity.  Therefore, $c \neq 0$.
Let $I$ be an ideal of $R$ that does not divide $c$ and consider
the congruence map
$$\eta: \Gamma \hookrightarrow \PSL(2,R) \rightarrow \PSL(2, R/I).$$
Then $\eta(\gamma) \notin \eta(A)$.

Let $H$ be a loxodromic subgroup of $\Gamma$. Let
$A$ be the maximal abelian subgroup of $\Gamma$ containing $H$,
and fix $\gamma \in \Gamma - H$.  By the case above, we may assume
that $\gamma \in A - H$.  Write $H = \langle h^m \rangle$, where
$h$ is a loxodromic element that generates $A$ and $m$ is a positive integer.
Since $\gamma \in A - H$, $\gamma = h^a$, for some integer $a$
that is not divisible by $m$.
By Corollary~\ref{orderm}, there exists a congruence map
$$\eta: \Gamma \hookrightarrow \PSL(2, R) \rightarrow \PSL(2, R/I)$$
such that the order of $\eta(h)$ is divisible by $m$.  Then
$\eta(\gamma) \notin \eta(H)$.

Let $H$ be a maximal cyclic parabolic subgroup of $\Gamma$.
Choose $k \in \Gamma$ such that $A = \langle h, k \rangle$ is the maximal
abelian subgroup of $\Gamma$ containing $H$, and fix $\gamma \in \Gamma - H$.
By the case above, we may assume that $\gamma \in A - H$.
Write $\gamma = h^ak^b$, $a, b \in \z$.
Since $\gamma \notin H$,  $b \neq 0$.
By Case 4 of Theorem~\ref{doublecoset}, for infinitely many primes $p$, there
exist congruence maps
$$\eta_p: \Gamma \hookrightarrow \PSL(2,R) \rightarrow \PSL(2, \Fp)$$
such that $\eta_p(A) \cong \z / p\z \oplus \z / p\z$.  If $p$
does not divide $b$, then $\eta_p(\gamma) \notin \eta_p(H)$.

Let $H$, $K$ and $g$ be given as in the statement of the corollary.
Since the congruence topology on $\Gamma$
is equivariant under left and right multiplication, it suffices to show
that the double coset $HK$ is closed in the congruence topology.
By the argument above, we may assume that $H$ and $K$ do not commute.
If $H$ and $K$ are both loxodromic, then the argument follows from Case 1
of Theorem~\ref{doublecoset}.  If $H$ is loxodromic and $K$
is a maximal parabolic subgroup of $\Gamma$, then the argument follows
from the proof of Lemma 3.2 of \cite{H2}.  If $H$ is loxodromic
and $K$ is a maximal cyclic parabolic subgroup of $\Gamma$, then
the argument follows from Case 4 of Theorem~\ref{doublecoset}.
\end{pf}


\section{Proof of conjugacy separability}\label{s: Conjugacy separability}

For a group $G$, we denote by $\widehat{G}$ the profinite
completion of $G$.  For a subgroup $H$ of $G$, we denote by
$\overline{H}$ the closure of $H$ in $\widehat{G}$. The proof of
Theorem \ref{t: Reducing cs to pieces of the torus decomposition}
relies on the main technical theorem of
\cite{wilton_profinite_2010}, which we state here for convenience.
We refer the reader to \cite{wilton_profinite_2010} for the
definitions of the terms \emph{efficient} and \emph{profinitely
acylindrical} used in it.

\begin{theorem}[\cite{wilton_profinite_2010}, Theorem 5.2]\label{t: Combination theorem}
Let $\G$ be a finite graph of groups with conjugacy separable vertex groups.
Let $G=\pi_1(\G)$, and suppose that the profinite topology on $G$ is
efficient and that $\G$ is profinitely 2-acylindrical.
For any vertex $v$ of $\G$ and incident edges $e$ and $f$, suppose
furthermore that the following conditions hold:
\begin{enumerate}
\item for any $g\in G_v$ the double coset $G_egG_f$ is separable in $G_v$;
\item the edge group $G_e$ is conjugacy distinguished in $G_v$;
\item the intersection of the closures of $G_e$ and $G_f$ in the profinite completion of $G_v$ is equal to
the profinite completion of their intersection.
That is, the natural map $\widehat{G_e\cap G_f}\to\overline{G}_e\cap\overline{G}_f$ is an isomorphism.
\end{enumerate}
Then $G$ is conjugacy separable.
\end{theorem}

\begin{definition}
A subgroup $H\subseteq G$ is called \emph{conjugacy distinguished} if,
whenever $g\in G$ is not conjugate into $H$, there is a finite quotient of
$G$ in which the image of $g$ is not conjugate into the image of $H$.
\end{definition}

As in \cite{wilton_profinite_2010}, we will apply Theorem \ref{t: Combination theorem} to the torus decomposition of $M$.  Henceforth, $M$ denotes a closed, orientable, Haken 3-manifold.
Let  $G=\pi_1(M)$ and let $\G$ be the graph-of-groups decomposition of
$G$ induced by the torus decomposition of $M$.

\begin{theorem}[\cite{wilton_profinite_2010}, Theorem A]\label{t: Efficiency}
For $M$ as above, the profinite topology of the fundamental group
of the graph of groups $\G$ is efficient.
\end{theorem}

It is convenient to make the extra assumption that every
Seifert-fibered piece of the torus decomposition of $M$ is
\emph{large}---that is, has a fundamental group that virtually
surjects a non-abelian free group. In this case, it turns out
that $\mathcal{G}$ is profinitely 2-acylindrical, and so Theorem
\ref{t: Combination theorem} applies.   In
\cite{wilton_profinite_2010}, the  remaining hypotheses of Theorem
\ref{t: Combination theorem} were checked for graph manifolds. We
will prove the corresponding results about the fundamental groups
of hyperbolic manifolds.

Therefore, we need to consider $G_v=\Gamma$, the fundamental group
of a finite-volume hyperbolic 3-manifold $N$. The incident edge
groups $G_e$ and $G_v$ are maximal parabolic subgroups of
$\Gamma$, which we shall denote $P$ and $Q$. The next lemma is a
consequence of Theorem~\ref{introdoublecoset}.

\begin{lemma}\label{l: Peripheral double-coset separability}
For any $g\in\Gamma$, the double coset $PgQ$ is separable in $\Gamma$.
\end{lemma}

\begin{lemma}\label{l: Hyp parabolic subgroups are conjugacy distinguished}
The subgroup $P$ is conjugacy distinguished in $\Gamma$.
\end{lemma}

\begin{pf}
As in the proof of Theorem \ref{doublecoset}, we may view $\Gamma \subset
\SL(2,R) \subset \SL(2,k)$, where $R$ is a finitely generated ring
contained in a number field $k$. Since $R$ is finitely generated,
$R \subset \okp$ for all but finitely many primes $\pp$ of $\okp$.
For each of these primes $\pp$ the residue class field map
$\eta_{\pp}: \okp \rightarrow F_{\pp}$ induces a group
homomorphism $$\overline{\eta}_{\pp}:\Gamma \hookrightarrow
\SL(2,\okp) \rightarrow \SL(2,F_{\pp}),$$ where $F_{\pp}$ is the
residue class field of $\okp$.

Let $\gamma$ be an element of $\Gamma$ that is not
conjugate into $P$.  If $\gamma$ is loxodromic, then
the square of the trace of $\gamma$, $\mathrm{tr}(\gamma)^2$,
is not equal to $4$.
Choose a prime $\pp$ of $\okp$ such that $R \subset \okp$ and
$\eta_{\pp}(\mathrm{tr}(\gamma)^2) \neq 4$.
Then $\overline{\eta}_{\pp}(\gamma) \notin
\overline{\eta}_{\pp}(P)$, as required.
Now suppose that $\gamma$ is parabolic.  As discussed in Case 3 of
Theorem~\ref{doublecoset},  we may view
$N$ as the interior of a compact manifold $N^{\prime}$
with a finite number of tori boundary components.  There is
a one to one correspondence between the boundary components of
$N^{\prime}$ and the conjugacy classes of maximal parabolic subgroups of
$\Gamma$.  Let $T_\gamma$ be the boundary component corresponding
to $\gamma$ and let $T_{P}$ be the boundary component corresponding to $P$.
Since $\gamma$ is not conjugate into $P$ and $P$ is a maximal parabolic subgroup,
$T_{\gamma} \neq T_{P}$.  By Thurston's Hyperbolic Dehn Surgery
Theorem, there exists a complete hyperbolic manifold $M$ of finite
volume obtained by Dehn surgery on $N^{\prime}$ along $T_{\gamma}$.
We may choose $M$ such that the image of $\gamma$ is non-trivial in $\pi_1(M)$.
Let $$\phi: \Gamma \cong \pi_1(N^{\prime}) \rightarrow \pi_1(M)$$ be the homomorphism
induced by inclusion $N^{\prime} \rightarrow M$.  Then $\phi(P)$ is a maximal parabolic
subgroup of $\pi_1(M)$ and $\phi(\gamma)$ is loxodromic.  The argument
then follows as above.
\end{pf}

\begin{lemma}\label{l: Distinguishing maps for the parabolic subgroups}
Let $P$ and $Q$ be non-conjugate maximal parabolic subgroups of
$\Gamma$.  There exists a sequence of group homomorphisms $f_n$
from $\Gamma$ to finite groups with the following properties.
\begin{enumerate}
\item For any finite index subgroup $K$ of $P$, there exists an $n$
with Ker$(f_n) \cap P \subset K$.
\item For every $n$, the intersection of $f_n(P)$ and $f_n(Q)$ is trivial.
\end{enumerate}
\end{lemma}

\begin{pf}
As in the proof of Lemma \ref{l: Hyp parabolic subgroups are conjugacy distinguished},
we may view $N$ as the interior of a compact manifold $N^{\prime}$ with a finite number of
tori boundary components, each corresponding to a conjugacy class of a maximal parabolic
subgroup of $\Gamma$. Let $T_P$  and $T_Q$ denote the boundary components of $N^{\prime}$
corresponding to $P$ and $Q$, respectively.  By assumption, $T_P \neq T_Q$.
Choose a basis $\{ p_1, p_2 \}$ of $P$ such that $N_{p_1}$ and $N_{p_2}$ are complete
hyperbolic $3$-manifolds of finite volume, where $N_{p_1}$ and $N_{p_2}$ are the manifolds
obtained by Dehn surgery on $N^{\prime}$ along $P$ sending $p_1$ and $p_2$, respectively,
to a meridian of the attached solid torus. Let $$\phi: \Gamma \cong \pi_1(N^{\prime}) \rightarrow \pi_1(N_{p_1})$$
be the homomorphism induced by inclusion.  Then $\phi(Q)$ is a maximal parabolic subgroup of $\pi_1(N_{p_1})$,
and $\phi(P)$ is a loxodromic subgroup of $\pi_1(N_{p_1})$ generated by $\phi(p_2)$.
We then proceed as in  Case 3 of Theorem~\ref{doublecoset}. We may view $\pi_1(N_{p_1}) \subset
\SL(2, R_1)$, where $R_1$ is a finitely  generated ring in a number field.  Moreover, we may assume that
$$\phi(p_2) =
\begin{pmatrix}
\lambda & 0 \\
0 & \lambda^{-1} \\ \end{pmatrix} \
\mathrm{and}$$
$$\phi(Q) =  \Big\{ \begin{pmatrix}
a & b \\ c & d \\ \end{pmatrix}
\begin{pmatrix}
1 & z \\ 0 & 1 \\ \end{pmatrix}
{\begin{pmatrix}
a & b \\ c & d \\ \end{pmatrix}}^{-1} \ \Big\vert \
z \in  A = \{ m + n \beta \ \vert \ m, n \in \z \} \Big\},$$
for some $$  \begin{pmatrix}
a & b \\ c & d \\ \end{pmatrix} \in \SL(2,R_1), \ \mathrm{and}
\ \lambda, \ \beta \in R_1.$$

Since $\pi_1(N_{p_1})$ is discrete, $\phi(P)$ and $\phi(Q)$ do not share a fixed point.
Therefore, $a$ and $c$ are non-zero. Fix a natural number $n$.  By Corollary~\ref{orderm},
there exist a finite field $F_n$ and a ring homomorphism $\eta_n: R_1 \rightarrow F_n$ such
that the multiplicative order of $\eta_n(\lambda)$ is divisible by $n$, $\eta_n(a) \neq 0$ and $\eta_n(c) \neq 0$.
This ring homomorphism induces a group homomorphism $$ \overline{\eta}_n: \pi_1(N_{p_1})
\hookrightarrow \SL(2, R_1) \rightarrow \SL(2,F_n)$$ such that the order of
$\overline{\eta}_n(\phi(p_2))$ is divisible by $n$.  Moreover, since $\eta_n(a) \neq 0$ and $\eta_n(c) \neq 0$,
the intersection of $\eta_n(\phi(\P))$ and $\eta_n(\phi(Q))$ is trivial.
Let $$g_n: \Gamma \cong \pi_1(N^{\prime}) \rightarrow \pi_1(N_{p_1}) \rightarrow \SL(2, F_n)$$
denote the composition $\overline{\eta}_n\circ\phi$.  Then the intersection of $g_n(P)$ and $g_n(Q)$ is trivial,
$g_n(p_1)$ is trivial, and $g_n(p_2)$ has order divisible by $n$.    By a similar argument, there exist a finite
field $L_n$ and a group homomorphism $$h_n: \Gamma \cong \pi_1(N^{\prime}) \rightarrow \pi_1(N_{p_2}) \rightarrow \SL(2, L_n)$$
such that the intersection of $h_n(P)$ and $h_n(Q)$ is trivial, $h_n(p_2)$ is trivial, and $h_n(p_1)$ has order divisible by $n$.
Let $K_n = \mathrm{Ker}(g_n) \cap \mathrm{Ker}(h_n)$, and let $$f_n: \Gamma \rightarrow \Gamma / K_n$$ denote the quotient map.
Then the intersection of $f_n(P)$ and $f_n(Q)$ is trivial, and Ker($f_n) \cap P = K_n \cap P \subset \langle np_1, np_2 \rangle = nP$.
The collection $\{ f_n \}$ satisfies the conditions above, since given a subgroup $K$ of finite index in $P$,
there exists a natural number $n$, such that $nP \subset K$.
\end{pf}

\begin{lemma}\label{l: Intersection of hyp closures}
Let $P,Q$ be distinct maximal parabolic subgroups of $\Gamma$. The intersection of the closures
$\overline{P}\cap\overline{Q}$ is trivial in the profinite completion $\widehat{\Gamma}$.
\end{lemma}

\begin{pf}
We first consider the case in which $P$ and $Q$ are not conjugate
in $\Gamma$.    By the universal property of the profinite
completion, the maps $f_n:\Gamma \to \Gamma/K_n$ in Lemma  \ref{l:
Distinguishing maps for the parabolic subgroups} extend to
continuous homomorphisms $\hat{f}_n$ from $\widehat{\Gamma}$. Let
$f=(\hat f_n):\widehat\Gamma \to \prod_n \Gamma/K_n$ be the
continuous homomorphism to the Cartesian product of $\Gamma/K_n$.
By Lemma  \ref{l: Distinguishing maps for the parabolic subgroups}
its restriction to $\overline P$ and $\overline Q$ is injective
and the image of $\overline P \cap\overline Q$ in $\Gamma/K_n$ is
trivial. Therefore $\overline P \cap\overline Q$ is trivial.

Suppose now that $P$ and $Q$ are conjugate in $\Gamma$.  As in Case 3 of the proof of Theorem \ref{doublecoset}, there is a subgroup $\Gamma_0$ of finite index in
$\Gamma$ such that $P_0=\Gamma_0\cap P$ and
$Q_0=\Gamma_0\cap Q$ are not conjugate in $\Gamma_0$.  By the
non-conjugate case, we have that $\overline{P}_0\cap
\overline{Q}_0=1$ in the profinite completion
$\widehat{\Gamma}_0$.  But $\widehat{\Gamma}_0$ is a subgroup of
finite index in $\widehat{\Gamma}$, and so
$\overline{P}\cap\overline{Q}$ is periodic in $\widehat{\Gamma}$.
But $\overline{P}\cong\widehat{\mathbb{Z}^2}$, which is
torsion-free, and so $\overline{P}\cap\overline{Q}=1$ as required.
\end{pf}

These lemmas complement the corresponding results for Seifert-fibered manifolds,
which we list below for convenience. As usual, by a \emph{peripheral}
subgroup of the fundamental group of a 3-manifold we mean a subgroup
that (up to conjugacy) corresponds to a boundary component of $N$.

\begin{lemma}[\cite{N}]\label{l: SF peripheral double coset separability}
Double cosets of peripheral subgroups in Seifert-fibered 3-manifold groups are separable.
\end{lemma}

\begin{lemma}[\cite{wilton_profinite_2010}, Lemma 5.3]
\label{l: SF parabolic subgroups are conjugacy distinguished}
Every peripheral subgroup of the fundamental group of a
Seifert-fibered 3-manifold group is conjugacy distinguished.
\end{lemma}

\begin{lemma}[\cite{wilton_profinite_2010}, Lemma 5.4]\label{l: Intersection of SF closures}
Let $N$ be a large Seifert-fibered 3-manifold and let $P,Q$
be distinct peripheral subgroups of $\pi_1(N)$.  Then $\widehat{P}\cap\widehat{Q}=\widehat{Z}$, the profinite
completion of the canonical normal cyclic subgroup of $\pi_1(N)$.
\end{lemma}

Proceeding exactly as in \cite{wilton_profinite_2010},
Lemmas \ref{l: Intersection of hyp closures} and \ref{l: Intersection of SF closures}
can be used together to generalize Lemma 5.5 of \cite{wilton_profinite_2010} as follows.

\begin{lemma}\label{l: Profinitely 2-acylindrical}
Let $M$ be a closed, orientable,
Haken 3-manifold in which every Seifert-fibered piece of the torus decomposition is large.
Let $G=\pi_1(M)$ and let $\mathcal{G}$ be the graph of groups induced by the torus decomposition of $M$.
Then $\mathcal{G}$ is profinitely 2-acylindrical.
\end{lemma}

When every Seifert-fibered piece of $M$ is large,
Theorem \ref{t: Reducing cs to pieces of the torus decomposition}
is a direct consequence of Theorem \ref{t: Combination theorem},
together with Theorem \ref{t: Efficiency} and Lemmas
\ref{l: Peripheral double-coset separability},
\ref{l: Hyp parabolic subgroups are conjugacy distinguished},
\ref{l: Intersection of hyp closures}, \ref{l: SF peripheral double coset separability},
\ref{l: SF parabolic subgroups are conjugacy distinguished},
\ref{l: Intersection of SF closures} and \ref{l: Profinitely 2-acylindrical}.
Of course, we also need the Geometrization Theorem (proved in the Haken case by Thurston and in full
by Perelman), which implies that every piece of the
torus decomposition is either Seifert-fibered or hyperbolic.

Finally, proceeding exactly as in the proof of Theorem D of
\cite{wilton_profinite_2010}, the full statement of
Theorem \ref{t: Reducing cs to pieces of the torus decomposition} follows.

\vskip.5in

\noindent
Department of Mathematics and Computer Science
\\ Emory University\\ Atlanta, GA 30322\\USA
\vskip.25in
\noindent
Department of Mathematics\\
University College London\\ Gower Street\\London\\WC1E 6BT\\UK
\vskip.25in
\noindent
Department of Mathematics\\
University of Brasilia\\
70910-900 Brasilia-DF\\
Brazil\\

\end{document}